\documentclass[12pt,twoside]{article}

\usepackage{epsfig,amsmath,amssymb}
\usepackage{amsfonts}
\usepackage{amssymb}
\usepackage{amssymb}
\usepackage{color}
\usepackage{enumerate}

\date{}

\begin{document}


\centerline{}

\centerline{}

\centerline {\Large{\bf New geometric constants of  
isosceles orthogonal type}}

\centerline{}

\centerline{\bf {}}

\centerline{Qi Liu,  Zhijian Yang, Yuankang Fu  and Yongjin Li$^*$}

\centerline{Department of Mathematics,
      Sun Yat-sen University}

\centerline{ Guangzhou, 510275, P. R. China}

\let\thefootnote\relax
\footnotetext{$^*$Corresponding author. E-mail: stslyj@mail.sysu.edu.cn}

\newtheorem{Theorem}{\quad Theorem}[section]

\newtheorem{Definition}[Theorem]{\quad Definition}

\newtheorem{Proposition}[Theorem]{\quad Proposition}

\newtheorem{Corollary}[Theorem]{\quad Corollary}

\newtheorem{Proof}[Theorem]{\quad Proof}

\newtheorem{Lemma}[Theorem]{\quad Lemma}

\newtheorem{Remark}[Theorem]{\quad Remark}

\newtheorem{Example}[Theorem]{\quad Example}

\newtheorem{Question}[Theorem]{\quad Question}

\begin{abstract}
Based on the parallelogram law and isosceles 
orthogonality, we define a new orthogonal geometric constant $ \Omega(X)$. 
We first discuss some basic properties of this new 
constant. Next,  we  consider the relation between the constant $\Omega(X)$ and the 
uniformly non-square property. 
Moreover, a generalized  constant $\Omega^{\prime}(X)$  is also introduced and some basic properties  are presented.
It is shown that, for a normed space, the constant 
value is equal to $1$ if and only if the norm  can be induced by the inner 
product. Finally, we  verify that this constant is closely related to   the well-known geometric constants through some inequalities.
\end{abstract}

{\bf Mathematics Subject Classification:} 46B20, 46C15.

{\bf Keywords:} Isosceles orthogonality, characterization
of Hilbert spaces, uniformly non-square, von Neumann-Jordan constant.

\section{Introduction and preliminaries}
\quad The characterization of inner product space has always been an active research area. Based on the parallelogram  law, scholars have obtained a series of characterization conditions, see \cite{FA2019,ERL1948,RRM1979}. 
The concept of orthogonality plays an important role in Euclidean geometry, which is closely related to many classical theorems. One of the themes of Banach space theory 
is to find suitable alternatives to this concept in Banach space.
Let $X$ be a real normed linear space and $x,y\in X$,  $x$ is said to be isosceles orthogonal to $y$ (denoted by $x\perp_Iy$) if
$\Vert x+y\Vert=\Vert x-y\Vert$. Birkhoff \cite{GB1935} introduced Birkhoff orthogonality: $x$ is said to be Birkhoff orthogonal to $y$  if  $\forall \in R, \Vert x+\lambda y\Vert \geq \Vert x\Vert$.
It is worth mentioning that the inner product space can also be described for elements that only satisfy some orthogonal conditions, see  \cite{AD1986,DSS2010,MB2006,SW2021}. This means that we do not need to consider all the elements in the space, but only those elements that meet the special orthogonal conditions. The following result is one of them, obtained by O.P. Kapoor et al.

\begin{Lemma} \label{Lemma}\cite{OPK}
Assume that $X$ be a normed space, $0<a,b<1$.
Then the following are equivalent: 

(i) $X$ be an inner product space. 

(ii) $x,y\in X$ and $
\|x+y\|^{2}+\|a x+b y\|^{2}=\|a x+y\|^{2}+\|x+b y\|^{2} 
$
implies $x\perp_Iy$.

(iii) $x,y\in X$ and $x\perp_Iy$ implies
$
\|x+y\|^{2}+\|ax+b y\|^{2}=\|a x+y\|^{2}+\|x+by\|^{2}
$.

Then $(i) \Rightarrow(i i) \Rightarrow(i i i)$, and $(i i i) \Rightarrow(i)$ when $a=b$.
\end{Lemma}

We now give some important definitions related to geometric constants.
Let real Banach space is represented by $X$ with $\operatorname{dim} X \geq 2$ throughout the paper. We will use $B_X$ and $S_X$ to denote the unit ball and unit sphere of $X$, respectively.

The  von Neumann-Jordan constant $C_{\mathrm{NJ}}(X)$  for  a Banach space $X$,  is defined by \cite{JAC1937}
$$C_{\mathrm{NJ}}(X)=\sup \left\{\frac{\|x+y\|^{2}+\|x-y\|^{2}}{2\left(\|x\|^{2}+\|y\|^{2}\right)}: x, y \in X, (x, y) \neq(0,0)\right\}.$$

We know that the  value of constant $C_{\mathrm{NJ}}(X)$ is $1$ if and only if  $X$ is a Hilbert space.
For more information on the properties of this classical and far-reaching constant, we recommend
\cite{KM2001,MH2021,TY1998}.

The James constant $J(X)$ of a Banach space $X$ is introduced by Gao and Lau \cite{GJ1990} as follows:
$$J(X)=\sup \{\min \{\Vert x+y\Vert,\Vert x-y\Vert \}:
x, y\in S_X\};$$
equivalently,
$$J(X)=\sup \{\min \{\Vert x+y\Vert,\Vert x-y\Vert \}:
x, y\in B_X\},$$
or 
$$
J(X)=\sup \left\{\|x+y\|: x, y \in S_{X}, x \perp_{I} y\right\}.
$$

Another non-square constant is Schäffer constants, which is denoted as
$$
S(X)=\inf \left\{\|x+y\|: x, y \in S_{X}, x \perp_{I} y\right\}.
$$
Noted that $J(X)S (X) = 2$  \cite{KM2001}. 

Recall  that the Banach space $X$ is called uniformly non-square
\cite{JRC1964}  if there exists a $\delta \in  (0, 1)$ such that for any
$x_1, x_2 \in S_X$ either $\frac{\Vert x_1+x_2\Vert}{2} \leq  1-\delta$ or $\frac{\Vert x_1-x_2\Vert}{2}\leq  1-\delta$. It is known,  that $X$ is uniformly non-square \cite{KM2001}
if and only if  $J(X)<2$.

The Clarkson modulus of convexity of a Banach space $X$  is the function
$\delta_X : [0, 2] \rightarrow [0, 1]$ defined by \cite{JAC1936}:
$$\delta_X(\varepsilon)=\inf \bigg\{1-\frac{\Vert x+y\Vert}{2}:x,y\in S_X, \Vert x-y\Vert\geq \varepsilon \bigg\}.$$

The constant $B R(X)$ is used to describe the difference  between   Birkhoff orthogonality 
and isosceles orthogonality
$$
B R(X)=\sup _{\alpha>0}\left\{\frac{\|x+\alpha y\|-\|x-\alpha y\|}{\alpha}: x, y \in S_{X}, x \perp_{B} y\right\}
$$
was studied by  P. L. Papini and S. Wu \cite{PLP2013}.

We collect some properties of this constant  for non-trivial Banach spaces
(see \cite{MH2018,PLP2013}):

(i) $0\leq BR(X)\leq 1$;

(ii) $BR(X)=0$ if and only if $X$ is an inner product space;

(iii) $
B R(X)=\sup \left\{\frac{\|x+y\|-\|x-y\|}{\|y\|}: x, y \in X, x, y \neq 0, x \perp_{B} y\right\}
$.

Ji et al
introduced the following geometric constant $D(X)$ as follows \cite{DJ2006}:
$$
D(X)=\inf \left\{\inf _{\lambda \in R}\|x+\lambda y\|: x, y \in S_{X}, x \perp_{I} y\right\},
$$
which is based on isosceles orthogonality.

\section{The Constant $\Omega(X)$}
\quad Combining the result of using the isosceles orthogonal condition to describe the inner product space in Lemma \ref{Lemma}, and considering the structure of some classical constants, we define the following orthogonal geometric constant. Compared with the classical constant, this new constant has both differences and similarities.

\begin{Definition}  For any Banach space $X$, we define
$$\begin{aligned}
      \Omega(X)=\sup\bigg\{\frac{\Vert x+2y\Vert^2+\Vert 2x+
y\Vert^2}{5\Vert x+y\Vert^2}:  x, y\in S_X,x\perp_Iy\bigg\}.
\end{aligned}$$
\end{Definition}

\begin{Proposition}
Let $X$ be a infinite-dimensional Banach space, then  $1\leq \Omega(X)\leq \frac{8}{5}$.
\end{Proposition}
{\bf Proof.}
According  Dvoretzki’s theorem, given
$\varepsilon>0$,  if the dimension of
$X$ is large enough, then there exists
a subspace $Y$ of  $X$, with $\operatorname{dim}(X)=2$, such that $|\Omega(X)-1|<\varepsilon$, this 
implies that $\Omega(X)\geq 1$ if $\operatorname{dim}(X)=\infty$.

For $x\perp_Iy$,  we can deduce that
$$\begin{aligned}
\frac{\|x+2 y\|^{2}+\|2 x+y\|^{2}}{5\|x+y\|^{2}} &\leq \frac{(\frac{3}{2}\Vert x+y\Vert+\frac{1}{2}\Vert x-y\Vert)^2+(\frac{3}{2}\Vert x+y\Vert+\frac{1}{2}\Vert x-y\Vert)^2}{5\|x+y\|^{2}}\\&=
 \frac{(2\Vert x+y\Vert)^2+(2\Vert x+y\Vert)^2}{5\|x+y\|^{2}}\\&=\frac{8}{5}.
\end{aligned}$$
$\hfill \blacksquare$

\begin{Example}
Let $X$ be $\mathbb{R}^2$  with the norm defined by
 $$\Vert (x_1,x_2)\Vert=\max \{|x_1|,|x_2|\}.$$
Letting $x_0=(1,0), y_0=(0,1)$, then $\Vert x_0+y_0\Vert=
\Vert x_0-y_0\Vert=1$, and $\|x_0+2 y_0\|=\|2 x_0+y_0\|=2$. Thus we can obtain that
$\Omega(X)=\frac{8}{5}$.
\end{Example}

\begin{Example}
      Let $X$ be $\mathbb{R}^{2}$ endowed with the $\ell_{\infty}-\ell_{1}$ norm
      $$
      \|x\|= \begin{cases}\|x\|_{1}, & x_{1} x_{2} \leq 0, \\ \|x\|_{\infty}, & x_{1} x_{2} \geq 0,\end{cases}
      $$
     then $\Omega(X)=\frac{49}{40}$.
      \end{Example}
      Assume that $x=\left(y_{1}, 1+y_{1}\right), y=\left(y_{2}, 1+y_{2}\right)$ with $-1 \leq y_{1} \leq y_{2} \leq 0$ and $x=\left(y_{1}, y_{1}-1\right), y=$ $\left(y_{2}, y_{2}-1\right)$ with $0 \leq y_{1} \leq y_{2} \leq 1$. The condition $\|x+y\|=\|x-y\|$ enforces that $\left|y_{1}-y_{2}\right|=2$, which contradicts.
      
      To estimate the constant value, we consider the following two cases.
      
      {\bf Case I}: Assume that $x=(x_{1},1), y=(1,y_{2})$ with $0\leq x_{1}\leq y_{2}\leq 1$.
    Since $\|x+y\|=\|x-y\|$, we have
      $$\|(x_{1}+1,y_{2}+1)\|=\|(x_{1}-1,1-y_{2})\|.$$
      Thus,
      $$1+y_{2}=(1-x_{1})+(1-y_{2}),$$
     and hence $x_{1}+2y_{2}=1$,
    which implies that  $y_{2}\in[\frac{1}{3},\frac{1}{2}]$.
Then we deduce that
      $$\begin{aligned}
      \frac{\|2x+y\|^2+\|x+2y\|^2}{5\Vert x+y\Vert^2}&=
      \frac{(2+y_{2})^2+(3-2y_{2})^2}{5(y_2+1)^2}\\
      &\leq\frac{(2+\frac{1}{3})^2+(3-\frac{2}{3})^2}{5(\frac{1}{3}+1)^2}
      \\&=\frac{49}{40},
      \end{aligned}$$
      and the maximum value it is attained  at the points  
     $x=(\frac{1}{3},1),y=(1,\frac{1}{3})$.
      
     {\bf Case II}: Assume that  $x=(x_{1},1),y=(y_{1},1+y_{1})$ with $-1\leq y_{1}\leq 0\leq x_{1}\leq 1$.
      Since  $\|x+y\|=\|x-y\|$, then
      $$\|(x_{1}+y_{1},2+y_{1})\|=\|(x_{1}-y_{1},-y_{1})\|.$$
      If $-x_{1}\leq y_{1}$, then we have
      $2+y_{1}=x_{1}-y_{1}$, 
     and hence $x_{1}-2y_{1}=2$, 
  we can obtain that $y_{1}\in[-\frac{2}{3},-\frac{1}{2}]$.
    Then we have
      $$\begin{aligned}
        \frac{\|2x+y\|^2+\|x+2y\|^2}{5\Vert x+y\Vert^2}&=\frac{\|(2x_{1}+y_{1},3+y_{1})\|^2+\|(x_{1}+2y_{1},3+2y_{1})\|^2}{5(x_1-y_1)^2}\\
      &=\frac{(3+y_{1})^2+(1-2y_{1})^2}{5(2+y_1)^2}\\
      &=\frac{49}{40}.
      \end{aligned}$$
      and the maximum value  it is attained only at the points  $x=(\frac{2}{3},1),y=(-\frac{2}{3},\frac{1}{3})$.
      
      If $y_{1}\leq-x_{1}$,  then we have
      $2-x_{1}=x_{1}-y_{1}$,
     and hence $2x_{1}-y_{1}=2$,
     which implies that $y_{1}\in[-1,-\frac{2}{3}]$.
     Thus,
  $$
        \frac{\|2x+y\|^2+\|x+2y\|^2}{5\Vert x+y\Vert^2}\leq \frac{49}{40}.
   $$
   and the maximum value  it is attained only at the points
    $x=(\frac{2}{3},1),y=(-\frac{2}{3},\frac{1}{3})$.
    
    To conclude, we have    $\Omega(X)=\frac{49}{40}$.

\section{Relations with other geometric constants and uniformly non-square property}
\begin{Theorem}\label{Theorem55}
      Let $X$ be a finite-dimensional  Banach space.  $X$ is not uniformly non-square if and only if  $\Omega(X)=\frac{8}{5}$.
      \end{Theorem}
      {\bf Proof.}
      Noted that, $X$ is not uniformly non-square and therefore
       there exist $x_n, y_n\in S_X$ for which
      $$
      \left\|x_n+y_n\right\| \rightarrow2,~\left\|x_n-y_n\right\| \rightarrow2~~(n\rightarrow\infty)
      .$$
      Letting  $u=\frac{x_n+y_n}{2},v=\frac{x_n-y_n}{2}$, then we have
      $\Vert u+v\Vert=\Vert u-v\Vert=1$.
      
      On the other hand, we have
      $$\begin{aligned}\Vert u+2v\Vert&= \Vert  \frac{3x_n-y_n}{2}\Vert \\&=\Vert 2x_n-2y_n-\frac{x_n-y_n}{2}+y_n\Vert\\&\geq \Vert 2x_n-2y_n\Vert-\Vert\frac{x_n-y_n}{2}\Vert-\Vert y_n\Vert\end{aligned}$$
      and $$\Vert 2x_n-2y_n\Vert -\Vert\frac{x_n-y_n}{2}\Vert-\Vert y_n\Vert\rightarrow2~~(n\rightarrow\infty).$$

      Combining with the fact that
      $$\Vert u+2v\Vert\leq \Vert u+v\Vert+\Vert v\Vert=2,$$
      we can deduce that  $\Vert u+2v\Vert\rightarrow2~~(n\rightarrow\infty)$. 
      Similarly, we have $\Vert 2u+v\Vert\rightarrow2~~(n\rightarrow\infty)$, which implies that  $\Omega(X)=\frac{8}{5}$.

      To prove the second part,  by $\Omega(X)=\frac{8}{5}$ and the fact that $X$ is finite dimensional we can  deduce that there exist
     $ x_{0}, y_{0} \in S_X$ with $x_{0} \perp_{I} y_{0}$ such that
      $$
      \frac{\left\|x_{0}+2 y_{0}\right\|^{2}+\left\|2 x_{0}+y_{0}\right\|^{2}}{5\left\|x_{0}+y_{0}\right\|^{2}}=\frac{8}{5}.
      $$
    Noted that $\left\|x_{0}+2 y_{0}\right\| \leq 2\left\|x_{0}+y_{0}\right\|,\left\|2 x_{0}+y_{0}\right\| \leq 2\left\|x_{0}+y_{0}\right\|$, and
     $$
     \frac{(2\left\|x_{0}+y_{0}\right\|)^{2}+(2\left\|x_{0}+y_{0}\right\|)^{2}}{5\left\|x_{0}+y_{0}\right\|^{2}} = \frac{8}{5}.
     $$
     Thus $\left\|x_{0}+2 y_{0}\right\|=2\left\|x_{0}+y_{0}\right\|$ and $\left\|2 x_{0}+y_{0}\right\|=2\left\|x_{0}+y_{0}\right\|$.

     Noted that
$$
\left\|x_{0}+2 y_{0}\right\| \leq\left\|y_{0}\right\|+\left\|x_{0}+y_{0}\right\|,\left\|2 x_{0}+y_{0}\right\| \leq\left\|x_{0}\right\|+\left\|x_{0}+y_{0}\right\| .
$$
Then we can obtain
$$
2\left\|x_{0}+y_{0}\right\| \leq\left\|y_{0}\right\|+\left\|x_{0}+y_{0}\right\|, 2\left\|x_{0}+y_{0}\right\| \leq\left\|x_{0}\right\|+\left\|x_{0}+y_{0}\right\|
$$
and hence
$$
\begin{aligned}
\max \left\{\left\|x_{0}+y_{0}\right\|,\left\|x_{0}-y_{0}\right\|\right\} &=\max \left\{\left\|x_{0}+y_{0}\right\|,\left\|x_{0}+y_{0}\right\|\right\} \\
& \leq \min \left\{\left\|x_{0}\right\|,\left\|y_{0}\right\|\right\}
\end{aligned}
$$

     Assume $X$ is uniformly non-square,  then there exists $\delta\in (0,2)$ such that
$$
\begin{array}{l}\min \left\{\left\|\frac{x+y}{\max \{\|x+y\|,\|x-y\|\}}+\frac{x-y}{\max \{\|x+y\|,\|x-y\|\}}\right\|,\right. \\ \left.~~~~~~~~~\left\|\frac{x+y}{\max \{\|x+y\|,\|x-y\|\}}-\frac{x-y}{\max \{\|x+y\|,\|x-y\|\}}\right\|\right\} \\ <2-\delta\end{array}
$$
for all $x,y\in X$, which implies that
$$
\min \left\{\left\|\frac{2x}{\max\{\Vert x+y\Vert,\Vert x-y\Vert \}}\right\|,\left\|\frac{2y}{\max\{\Vert x+y\Vert,\Vert x-y\Vert \}}\right\|\right\}<2-\delta.
$$
Thus 
$$
\min \{\|x\|,\|y\|\}<(1-\frac{\delta}{2}) \max\{\|x+y\|,\|x-y\|\}
$$
and hence 
$$
 \max\{\|x+y\|,\|x-y\|\}> \frac{2}{2-\delta}\min \{\|x\|,\|y\|\}.
$$This contradicts the fact  
$$
\max\{\|x_0+ y_0\|,\|x_0- y_0\|\}\leq \min\{\Vert x_0\Vert, \Vert y_0\Vert\},
$$   and 
thus, we complete the proof.
$\hfill \blacksquare$

      \begin{Proposition} Let $X$ be a non-trivial Banach space, then
$$
\frac{8}{5}+\frac{2}{5}\frac{1}{J(X)^2}-\frac{8}{5}\frac{1}{J(X)}\leq \Omega(X) \leq \frac{2}{5}+\frac{J(X)^{2}}{10}+\frac{2 J(X)}{5}.
$$
\end{Proposition}
{\bf
Proof.} For any $x, y \in S_{X}, x \perp_{I} y$, we have
$$
\begin{aligned}
\frac{2 \min \{\|x+y\|,\|x-y\|\}}{\Vert x+y\Vert} &=\frac{\min \{\|2 x+2 y\|,\|2 x+2 y\|\}}{\Vert x+y\Vert} \\
& \leq \frac{\min \{1+\|x+2 y\|, 1+\|2 x+y\|\}}{\Vert x+y\Vert}
 \\
&=\frac{1+\min \{\|x+2 y\|,\|2 x+y\|\}}{\Vert x+y\Vert} \\
& \leq \frac{1+\sqrt{\|x+2 y\| \cdot\|2 x+y\|}}{\Vert x+y\Vert}
\\
& \leq \frac{1}{\Vert x+y\Vert}+\sqrt{\frac{\|x+2 y\|^2+\|2 x+y\|^2}{2\Vert x+y\Vert^2}},
\end{aligned}
$$
which implies that
$$\sqrt{\frac{2}{5}}\cdot\frac{2 \min \{\|x+y\|,\|x-y\|\}-1}{\Vert x+y\Vert}\leq \sqrt{\frac{\|x+2 y\|^2+\|2 x+y\|^2}{5\Vert x+y\Vert^2}}$$
and hence
$$\sqrt{\frac{8}{5}}-\sqrt{\frac{2}{5}}\frac{1}{J(X)}\leq \sqrt{\Omega(X)}.$$
Then we can obtain that
$$\frac{8}{5}+\frac{2}{5}\frac{1}{J(X)^2}-\frac{8}{5}\frac{1}{J(X)}\leq \Omega(X).$$

On the other hand, for any $x, y \in S_{X}, x \perp_{I} y$,  we can obtain that
$$
\begin{aligned}
\frac{\|2 x+y\|^2+\|x+2 y\|^2}{5\Vert x+y\Vert^2} & \leq\frac{(\|x+y\|+\|x\|)^2+(\|x+y\|+\|y\|)^2}{5\Vert x+y\Vert^2} \\
& =\frac{2}{5}+\frac{4}{5\Vert x+y\Vert}+\frac{2}{5\Vert x+y\Vert^2}\\
& \leq\frac{2}{5}+\frac{4}{5S(X)}+\frac{2}{5S(X)^2}.
\end{aligned}
$$
According $S(X) \cdot J(X)=2$, we can obtain that
$$
\frac{\|2 x+y\|^{2}+\|x+2 y\|^{2}}{5\|x+y\|^{2}} \leq \frac{2}{5}+\frac{J(X)^{2}}{10}+\frac{2 J(X)}{5},
$$
as desired.

\section{Equivalent form of $\Omega(X)$ in symmetric Minkowski plane}
Let $X$ be a Minkowski plane, if there exist $e_{1}, e_{2} \in S_X$ such that:
$$
\left\|e_{1}+t e_{2}\right\|=\left\|e_{1}-t e_{2}\right\|=\left\|e_{2}+t e_{1}\right\|=\left\|e_{2}-t e_{1}\right\|
$$
holds for all $t \in \mathbb{R}$, then we call $X$ a symmetric Minkowski plane \cite[p.5]{JD} and $\left\{e_{1}, e_{2}\right\}$ a pair of axes of $X$.

By following the ideas in \cite{JD,YG2007}, 
we will give the equivalent form of $\Omega(X)$ in symmetric Minkowski plane. 

\begin{Lemma}\cite{JD}\label{L1}
Let $X$ be a symmetric Minkowski plane, $\{e_1, e_2\}$ be a pair of axes of $X$. Then $\forall x,y \in S_X, x=\alpha e_1+\beta e_2, x\perp_I y$ iff $y=\pm(-\beta e_1+\alpha e_2)$.
\end{Lemma}

\begin{Proposition}
 Let $X$ be a symmetric Minkowski plane, ${e_1, e_2}$ be a pair of axes of $X$. Then
 $$\Omega(X)=\max\left\{\frac{f(t)^2+f(-t)^2}{5g(t)^2}~:~0\leq t<\infty\right\},$$
 where
 $$f(t)=\|(1+2t)e_1+(2-t)e_2\|,$$
 $$g(t)=\|(1+t)e_1+(1-t)e_2\|.$$
\end{Proposition}

{\bf Proof}
{\bf Step 1}: For any  $x=\alpha e_1+\beta e_2\in S_X$,   we can suppose that $\alpha\neq0$. Then, for any $x=\alpha e_1+\beta e_2\in S_X$, we have
$$\left \|e_1+\frac{\beta}{\alpha} e_2\right\|=\left\|\frac{1}{\alpha}x\right\|=\frac{1}{|\alpha|},$$
which follows that
\begin{equation}\label{E1}
  x=\frac{e_1+\frac{\beta}{\alpha} e_2}{\|e_1+\frac{\beta}{\alpha} e_2\|}sgn \alpha.
\end{equation}
By   Lemma \ref{L1}, we have
\begin{equation}\label{E2}
 y=\pm\frac{\frac{\beta}{\alpha}e_1- e_2}{\|\frac{\beta}{\alpha}e_1-e_2\|},
\end{equation}
where $y\in S_X$ and $x\perp_I y$. For the convenience, we denote $\frac{\frac{\beta}{\alpha}e_1- e_2}{\|\frac{\beta}{\alpha}e_1-e_2\|}$ and $-\frac{\frac{\beta}{\alpha}e_1- e_2}{\|\frac{\beta}{\alpha}e_1-e_2\|}$ by $y_x, \overline{y}_x$, respectively.

{\bf Step 2}: Now,  for any  $x=\alpha e_1+\beta e_2\in S_X$, according (\ref{E1}) and (\ref{E2}), we have
\begin{align*}
 \|x+2y_x\|&=\left\|\frac{e_1+\frac{\beta}{\alpha} e_2}{\|e_1+\frac{\beta}{\alpha} e_2\|}sgn \alpha+2\frac{\frac{\beta}{\alpha}e_1- e_2}{\|\frac{\beta}{\alpha}e_1-e_2\|}\right\| \\
  & =\left\|\frac{e_1+\frac{\beta}{\alpha} e_2}{\|e_1+\frac{\beta}{\alpha} e_2\|}+2\frac{\frac{\beta}{\alpha}e_1- e_2}{\|e_1+\frac{\beta}{\alpha}e_2\|}\right\| \\
  &=\frac{\|(1+2\frac{\beta}{\alpha})e_1+(\frac{\beta}{\alpha}-2)e_2\|}{\|e_1+\frac{\beta}{\alpha}e_2\|}\\
    &=\frac{\|(1+2\frac{\beta}{\alpha})e_1+(2-\frac{\beta}{\alpha})e_2\|}{\|e_1+\frac{\beta}{\alpha}e_2\|},
\end{align*}
and,
\begin{align*}
 \|2x+\overline{y}_x\|&=\left\|2\frac{e_1+\frac{\beta}{\alpha} e_2}{\|e_1+\frac{\beta}{\alpha} e_2\|}sgn \alpha-\frac{\frac{\beta}{\alpha}e_1- e_2}{\|\frac{\beta}{\alpha}e_1-e_2\|}\right\| \\
  & =\left\|\frac{2e_1+2\frac{\beta}{\alpha} e_2}{\|e_1+\frac{\beta}{\alpha} e_2\|}-\frac{\frac{\beta}{\alpha}e_1- e_2}{\|e_1+\frac{\beta}{\alpha}e_2\|}\right\| \\
  &=\frac{\|(2-\frac{\beta}{\alpha})e_1+(2\frac{\beta}{\alpha}+1)e_2\|}{\|e_1+\frac{\beta}{\alpha}e_2\|}\\
    &=\frac{\|(1+2\frac{\beta}{\alpha})e_1+(2-\frac{\beta}{\alpha})e_2\|}{\|e_1+\frac{\beta}{\alpha}e_2\|},
\end{align*}
which means that
\begin{equation}\label{E3}
\|x+2y_x\|=\|2x+\overline{y}_x\|=\frac{\|(1+2\frac{\beta}{\alpha})e_1+(2-\frac{\beta}{\alpha})e_2\|}{\|e_1+\frac{\beta}{\alpha}e_2\|}.
\end{equation}
Similarly, we can get
\begin{equation}\label{E4}
\|x+2\overline{y}_x\|=\|2x+y_x\|=\frac{\|(1-2\frac{\beta}{\alpha})e_1+(2+\frac{\beta}{\alpha})e_2\|}{\|e_1+\frac{\beta}{\alpha}e_2\|}.
\end{equation}
According (\ref{E1}) and (\ref{E2}), we also have
\begin{align*}
\|x+y_x\|&=\left\|\frac{e_1+\frac{\beta}{\alpha} e_2}{\|e_1+\frac{\beta}{\alpha} e_2\|}sgn \alpha+\frac{\frac{\beta}{\alpha}e_1- e_2}{\|\frac{\beta}{\alpha}e_1-e_2\|}\right\| \\
 & =\left\|\frac{e_1+\frac{\beta}{\alpha} e_2}{\|e_1+\frac{\beta}{\alpha} e_2\|}+\frac{\frac{\beta}{\alpha}e_1- e_2}{\|e_1+\frac{\beta}{\alpha}e_2\|}\right\| \\
  &=\frac{\|(1+\frac{\beta}{\alpha})e_1+(\frac{\beta}{\alpha}-1)e_2\|}{\|e_1+\frac{\beta}{\alpha}e_2\|}\\
 &=\frac{\|(1+\frac{\beta}{\alpha})e_1+(1-\frac{\beta}{\alpha})e_2\|}{\|e_1+\frac{\beta}{\alpha}e_2\|}.
\end{align*}
Further, since $x\perp_I y_x$, we have
\begin{equation}\label{E5}
  \|x+\overline{y}_x\|=\|x-y_x\|=\|x+y_x\|=\frac{\|(1+\frac{\beta}{\alpha})e_1+(1-\frac{\beta}{\alpha})e_2\|}{\|e_1+\frac{\beta}{\alpha}e_2\|}.
\end{equation}

{\bf Step 3}: Further, according (\ref{E3}), (\ref{E4}) and (\ref{E5}), we have
 \begin{align}\label{E6}
 &\left\{\frac{\|x+2y_x\|^2+\|2x+y_x\|^2}{5\|x+y_x\|^2}~:~x\in S_X\right\}\\ \nonumber =&\left\{\frac{(\frac{\|(1+2\frac{\beta}{\alpha})e_1+(2-\frac{\beta}{\alpha})e_2\|}{\|e_1+\frac{\beta}{\alpha}e_2\|})^2+(\frac{\|(1-2\frac{\beta}{\alpha})e_1+(2+\frac{\beta}{\alpha})e_2\|}{\|e_1+\frac{\beta}{\alpha}e_2\|})^2}{5(\frac{\|(1+\frac{\beta}{\alpha})e_1+(1-\frac{\beta}{\alpha})e_2\|}{\|e_1+\frac{\beta}{\alpha}e_2\|})^2}~:~
 -\infty<\frac{\beta}{\alpha}<\infty\right\}\\\nonumber
 =& \left\{\frac{(\frac{\|(1-2\frac{\beta}{\alpha})e_1+(2+\frac{\beta}{\alpha})e_2\|}{\|e_1+\frac{\beta}{\alpha}e_2\|})^2+(\frac{\|(1+2\frac{\beta}{\alpha})e_1+(2-\frac{\beta}{\alpha})e_2\|}{\|e_1+\frac{\beta}{\alpha}e_2\|})^2}{5(\frac{\|(1+\frac{\beta}{\alpha})e_1+(1-\frac{\beta}{\alpha})e_2\|}{\|e_1+\frac{\beta}{\alpha}e_2\|})^2}~:~
 -\infty<\frac{\beta}{\alpha}<\infty\right\}\\\nonumber
 =&\left\{\frac{\|x+2\overline{y}_x\|^2+\|2x+\overline{y}_x\|^2}{5\|x+\overline{y}_x\|^2}~:~x\in S_X\right\}.
 \end{align}

{\bf Step 4}: Since $X$ is finite-dimensional space, it is clearly that 
$$\begin{aligned} \Omega(X)&=\sup\left\{\frac{\|x+2y\|^2+\|2x+y\|^2}{5\|x+y\|^2}~:~x,y\in S_X, x\perp_I y\right\}\\&=\max\left\{\frac{\|x+2y\|^2+\|2x+y\|^2}{5\|x+y\|^2}~:~x,y\in S_X, x\perp_I y\right\}.
	\end{aligned}$$
	From (\ref{E6}), we have
	\begin{align*}
		\Omega(X)
		&=\max\left\{\frac{\|x+2y\|^2+\|2x+y\|^2}{5\|x+y\|^2}~:~x,y\in S_X, x\perp_I y\right\} \\
		& =\max\left\{\left\{\frac{\|x+2y_x\|^2+\|2x+y_x\|^2}{5\|x+y_x\|^2}~:~x\in S_X\right\} \bigcup \left\{\frac{\|x+2\overline{y}_x\|^2+\|2x+\overline{y}_x\|^2}{5\|x+\overline{y}_x\|^2}~:~x\in S_X\right\}\right\}\\
		&=\max\left\{\frac{\|x+2y_x\|^2+\|2x+y_x\|^2}{5\|x+y_x\|^2}~:~x\in S_X\right\}\\
		&=\max\left\{\frac{\bigg(\frac{\|(1+2\frac{\beta}{\alpha})e_1+(2-\frac{\beta}{\alpha})e_2\|}{\|e_1+\frac{\beta}{\alpha}e_2\|}\bigg)^2+\bigg(\frac{\|(1-2\frac{\beta}{\alpha})e_1+(2+\frac{\beta}{\alpha})e_2\|}{\|e_1+\frac{\beta}{\alpha}e_2\|}\bigg)^2}{5\bigg(\frac{\|(1+\frac{\beta}{\alpha})e_1+(1-\frac{\beta}{\alpha})e_2\|}{\|e_1+\frac{\beta}{\alpha}e_2\|}\bigg)^2}~:~
		-\infty<\frac{\beta}{\alpha}<\infty\right\}.
	\end{align*}

For the converence,  letting $t=\frac{\beta}{\alpha}$,
we can get
	\begin{align*}
	\Omega(X)
	&=\max\left\{\frac{\|x+2y\|^2+\|2x+y\|^2}{5\|x+y\|^2}~:~x,y\in S_X, x\perp_I y\right\}\\
	&=\max\left\{\frac{(\|(1+2t)e_1+(2-t)e_2\|)^2+(\|(1-2t)e_1+(2+t)e_2\|)^2}{5(\|(1+t)e_1+(1-t)e_2\|)^2}~:~
	-\infty<t<\infty\right\}\\
	&=\max\left\{\frac{f(t)^2+f(-t)^2}{5g(t)^2}~:~-\infty<t<\infty\right\}.
\end{align*}

On the other hand, we have
	$$\frac{f(t)^2+f(-t)^2}{5g(t)^2}=\frac{f(-t)^2+f(t)^2}{5g(-t)^2},$$
which implies that
	$$ \Omega(X)=\max\left\{\frac{f(t)^2+f(-t)^2}{5g(t)^2}~:~0\leqslant t<\infty\right\}. $$

\section{The constant $\Omega^{\prime}(X)$}
Small changes in the range of geometric constants affect the constants. Without considering the conditions $x,y \in S_X$ in the definition of constant $\Omega(X)$, the following more general definition of constant can be obtained.
\begin{Definition}  For any Banach space $X$, we define
      $$\begin{aligned}
            \Omega^{\prime}(X)=\sup\bigg\{\frac{\Vert x+2y\Vert^2+\Vert 2x+
      y\Vert^2}{5\Vert x+y\Vert^2}:  x\perp_Iy\bigg\}.
      \end{aligned}$$
      \end{Definition}
      \begin{Remark}
       We can also view  $\Omega^{\prime}(X)$ as following:
      $$\begin{aligned}
            \Omega^{\prime}(X)&=\sup\bigg\{\frac{\Vert x+2y\Vert^2+\Vert 2x+
      y\Vert^2}{5\Vert x+y\Vert^2}: x\in S_X, y\in B_X, x\perp_Iy\bigg\}
      \\&=\sup\bigg\{\frac{2}{5}\frac{\Vert x+y\Vert^2+\Vert 2x+
      y\Vert^2}{\Vert x+y\Vert^2+\Vert x-y\Vert^2}: x\in S_X, y\in B_X, x\perp_Iy\bigg\}.
      \end{aligned}$$
      \end{Remark}

      \begin{Proposition}
      Let $X$ be a Banach space, then  $1\leq\Omega^{\prime}(X)\leq \frac{8}{5}$.
      \end{Proposition}
      {\bf Proof.} Letting $y=0$, then $x\perp_Iy$ and  hence
      $$\Omega^{\prime}(X)\geq \frac{\Vert x\Vert^2+\Vert 
      2x\Vert^2}{5\Vert x\Vert^2}=1.$$

      For $x\perp_Iy$,  we can deduce that
      $$\begin{aligned}
      \frac{\|x+2 y\|^{2}+\|2 x+y\|^{2}}{5\|x+y\|^{2}} &\leq \frac{(\frac{3}{2}\Vert x+y\Vert+\frac{1}{2}\Vert x-y\Vert)^2+(\frac{3}{2}\Vert x+y\Vert+\frac{1}{2}\Vert x-y\Vert)^2}{5\|x+y\|^{2}}\\&=
       \frac{(2\Vert x+y\Vert)^2+(2\Vert x+y\Vert)^2}{5\|x+y\|^{2}}\\&=\frac{8}{5}.
      \end{aligned}$$
      
      Before continuing to prove that the constant $\Omega^{\prime}(X)$ can describe the inner product space, the following  technical lemma is  established in the 
      proofs of \cite{BC1984}, is required.
      \begin{Lemma}\label{Lemma2}
      \cite{BC1984} Let $X$ be a real normed linear space. Then 
       $\Vert \cdot\Vert$ derives from an inner product if and 
       only if for all $x,y$ in $S_X$ there exist 
       $\alpha,\beta \neq  0$ such that
       $$
      \|\alpha x+\beta y\|^{2}+\|\alpha x-\beta y\|^{2} \sim 2\left(\alpha^{2}+\beta^{2}\right)
      ,$$
      where $\sim$  stands for $=$, $\leq $ or $\geq$.
      \end{Lemma}
      
      \begin{Proposition}
      Let $(X,\Vert \cdot\Vert)$ be a normed space. Then, $\Vert \cdot\Vert $ 
       derives from an inner product if and only if
       $\Omega^{\prime}(X)=1$.
      \end{Proposition}
      
      {\bf Proof} Noted that, since $\Omega^{\prime}(X)=1$, we can obtain that 
      $$\frac{\Vert x+2y\Vert^2+\Vert 2x+
      y\Vert^2}{5\Vert x+y\Vert^2}\leq 1$$
      for all $\Vert x\Vert=\Vert y\Vert$ with $x\perp_Iy$.
      
      Letting $x,y\in S_X$, then  $x+y\perp_Ix-y$.  We can deduce that
      $$\begin{aligned}
      \|2(x+y)+(x-y)\|^{2}+\|(x+y)+2(x-y)\|^{2}\leq5\|x+y+x-y\|^{2}
      \end{aligned}$$
      for all $x,y\in S_X$,
      which implies that 
      $$\begin{aligned}
      \|3x+y\|^{2}+\|3x-y\|^{2}\leq20
      \end{aligned}$$
      for all $x,y\in S_X$.
       Setting $\alpha=3,\beta=1$ in 
      Lemma \ref{Lemma2}, we can deduce that  $\Vert \cdot\Vert$ derives from an inner product.
      
      The rest of the proof, we
      assume that $a=\frac{1}{2}, b=\frac{1}{2}$ in Lemma \ref{Lemma} (iii), then we get
      $$
      \|x+y\|=\|x-y\| \Rightarrow   4\|x+y\|^{2}+\| x+ y\|^{2}=\| x+2y\|^{2}+\|2x+ y\|^{2} 
      ,$$
      which implies that $\Omega^{\prime}(X)=1$.
      $\hfill \blacksquare$

    \begin{Theorem}
Let $X$ be a Banach space. Then,
$$\Omega^{\prime}(X)\leq \frac{2}{5}C_{NJ}(X)+\frac{4}{5}.
$$
\end{Theorem}
{\bf Proof.} 
Noted that, $C_{N J}(X)$ can be written in the following equivalent form
$$
C_{N J}(X)=\sup \bigg\{ \frac{2\left(\|x\|^{2}+\|y\|^{2}\right)}{\|x+y\|^{2}+\|x-y\|^{2}}
:(x,y)\neq (0,0)\bigg\}.
$$

For $x\perp_Iy$, by applying the triangle inequality, we can obtain the following estimate:
$$\begin{aligned}
&\frac{\|x+2 y\|^{2}+\|2 x+y\|^{2}}{5\|x+y\|^{2}}\\&=2\cdot
\frac{\|x+2 y\|^{2}+\|2 x+y\|^{2}}{5\|x+y\|^{2}+5\|x-y\|^{2}}
\\&\leq 2\cdot\frac{(\Vert x+y\Vert+\Vert y\Vert)^2+
(\Vert x+y\Vert+\Vert x\Vert)^2}{5\|x+y\|^{2}+5\|x-y\|^{2}}
\\&\leq 2\cdot\frac{2(\Vert x+y\Vert^2+\Vert x-y\Vert^2)+2(\Vert x\Vert^2+\Vert y\Vert^2)}{5\|x+y\|^{2}+5\|x-y\|^{2}}
\\&\leq \frac{2}{5}C_{NJ}(X)+\frac{4}{5},
\end{aligned}$$
as desired.
$\hfill \blacksquare$

Next, we will introduce  constant $\gamma_{X} (t)$, which not only plays an important role in estimating  the von Neumann-Jordan constant, but is also 
closely related to  constant $\Omega^{\prime}(X)$. For more detail, we refer  the reader
to the papers \cite{MH2015,YC2006}.
\begin{Definition}\cite{YC2006}
Let $X$ be a Banach space. The function 
$\gamma_{X} (t ) : [0, 1] \rightarrow [1, 4]$ is defined by

$$\gamma_{X}(t)=\sup\left\{\frac{\|x+ty\|^{2}+\|x-ty\|^{2}}{2}:x,y\in S_{X}\right\}.$$
\end{Definition}
\begin{Proposition}\label{Proposition11}
Let $X$ be a Banach space. Then
$\Omega^{\prime}(X)=\frac{9}{10}\gamma_{X}(\frac{1}{3})$.
\end{Proposition}
{\bf Proof.} We divide the proof into two steps:

{\bf Step~1} For any  $x\perp_{I} y$, let $u=\frac{x+y}{2} \text{ and } v=\frac{x-y}{2}$.  We can obtain that  $$\begin{aligned}\frac{\|x+2y\|^{2}+\|2x+y\|^{2}}{5\|x+y\|^{2}}&=\frac{\|3u-v\|^{2}+\|3u+v\|^{2}}{20\|u\|^{2}}\\&=\frac{9}{20}\frac{\|u-\frac{1}{3}v\|^{2}+\|u+\frac{1}{3}v\|^{2}}{\|u\|^{2}},
\end{aligned}$$
where $\|u\|=\|v\|$.
 
 Letting $x'=\frac{u}{\|u\|}, y'=\frac{v}{\|v\|}$,  we have 
$$\frac{\|u-\frac{1}{3}v\|^{2}+\|u+\frac{1}{3}v\|^{2}}{\|u\|^{2}}=\|x'-\frac{1}{3} y'\|^{2}+\|x'+\frac{1}{3} y'\|^{2}\leq 2\gamma_{X}(\frac{1}{3}),$$
which implies that $\Omega^{\prime}(X)\leq \frac{9}{10}\gamma_{X}(\frac{1}{3}).$

{\bf Step~2} For any $x,y\in S_{X}$, let $u=\frac{x+y}{2}$ and $v=\frac{x-y}{2}$. 
It is clear that
$\|u+v\|=\|u-v\|=1$.
Then we can deduce that  
$$\begin{aligned}\frac{\|x-\frac{1}{3}y\|^{2}+\|x+\frac{1}{3}y\|^{2}}{2}&=\frac{\|u+v-\frac{1}{3}(u-v)\|^{2}+\|u+v+\frac{1}{3}(u-v)\|^{2}}{2\|u+v\|^{2}}
\\&=\frac{2}{9}\cdot\frac{\|u+2v\|^{2}+\|2u+v\|^{2}}{\|u+v\|^2}\\&\leq \frac{10}{9}\Omega^{\prime}(X),\end{aligned}$$
and hence $\frac{9}{10}\gamma_{X}(\frac{1}{3})\leq \Omega^{\prime}(X)$.
$\hfill \blacksquare$

 \section*{Data Availability Statement}
All type of data used for supporting the conclusions of this article is included in the article and also is cited at relevant places within the text as references.

 \section*{Conflict of interest}
 The authors declare that they have no conflict of interest.

 \section*{ Funding Statement}
This work was supported by the National Natural Science Foundation 
of P. R. China (Nos. 11971493 and 12071491).

\end{document}